\documentclass{amsart}[11pt]

\usepackage{amsmath, dsfont, amsthm, amssymb, latexsym, graphicx, fullpage}

\def\identity{\mathds 1}
\def\CC{\mathbb C}
\def\bra{\left\langle}
\def\ket{\right\rangle}
\DeclareMathOperator{\rank}{\bf rk}
\DeclareMathOperator{\supp}{\bf supp}
\DeclareMathOperator{\trace}{\bf tr}
\newcommand{\op}{\operatorname}

\numberwithin{equation}{section}

\newtheorem{intro-thm}{Theorem}
\newtheorem{intro-cor}{Corollary}

\newtheorem*{thm-restate-main}{Theorem~\ref{thm:main}}
\newtheorem*{thm-restate-func}{Theorem~\ref{thm:func}}
\newtheorem*{cor-restate-cor1}{Corollary~\ref{cor:cor1}}
\newtheorem*{cor-restate-cor2}{Corollary~\ref{cor:cor2}}
\newtheorem{thm}{Theorem}[subsection]

\newtheorem{prop}{Proposition}[subsection]

\begin{document}

\title{Uncertainty Principles for Compact Groups}

\author{Gorjan Alagic} 
\address{Department of Mathematics, University of Connecticut} 
\email{alagic@math.uconn.edu} 
\author{Alexander Russell} 
\address{Department of Computer Science and Engineering, University of Connecticut} 
\email{acr@cse.uconn.edu} 

\thanks{Alexander Russell gratefully acknowledges support from NSF grants CCR-0220070, EIA-0218563, and CCF-0524613, and ARO contract W911NF-04-R-009.}


\begin{abstract}
We establish an operator-theoretic uncertainty principle over arbitrary compact groups, generalizing several previous results. As a consequence, we show that every nonzero square-integrable function $f$ on a compact group $G$ satisfies $\mu(\supp f)\cdot \sum_{\rho \in \hat G} d_\rho \rank \hat f(\rho) \geq 1$. For finite groups, our principle implies the following: if $\op P$ and $\op R$ are projection operators on the group algebra $\CC G$ such that $\op P$ commutes with projection onto each group element, and $\op R$ commutes with left multiplication, then $\| \op{PR} \|^2 \leq (\rank \op P \rank \op R) / |G|$.
\end{abstract}

\maketitle

\section{Introduction}
Uncertainty principles assert, roughly, that \emph{a function and its Fourier transform cannot simultaneously be highly concentrated}. One example of this is the well-known Heisenberg uncertainty principle concerning position and momentum wavefunctions in quantum physics. Several uncertainty principles have been formulated for complex-valued functions on groups. For finite abelian groups, perhaps the most basic of these is an inequality which relates the sizes of the supports of $f$ and its transform $\hat f$ to the size of the group. It states that
\begin{equation}\label{finite-func}
|\supp f|| \supp \hat f| \geq |G|\,,
\end{equation}
unless $f$ is identically zero \cite{Terras}. One can generalize this inequality by establishing an analogous statement concerning the associated projection operators on the group algebra $\CC G$. Specifically, let $\op P$ be the operator which, when expressed in the group basis, projects to the support of $f$; let $\op R$ be the operator which, when expressed in the Fourier basis, projects to the support of $\hat f$. Since the operator norm of $\op {PR}$ is equal to $1$, we can rewrite \eqref{finite-func} as
$\| \op {PR}\|^2 \leq |\supp  f||\supp  \hat{f}|/|G|$.
It is natural to ask if a similar fact holds for \emph{any} pair $\op P_S$ and $\op R_T$, where $\op P_S$ projects (in the group basis) to $S \subset G$ and $\op R_T$ projects (in the Fourier basis) to $T \subset \widehat G.$ Indeed, a generalized uncertainty principle for finite abelian groups states that
\begin{equation}\label{finite-operator}
\|\op P_S\op R_T\|^2 \leq  \frac{|S||T|}{|G|}\,.
\end{equation}
The above principles were proved for the real and finite cyclic case by Donoho and Stark~\cite{Donoho} and extended to locally compact abelian groups by Smith~\cite{Smith}. See also Terras's presentation~\cite{Terras} of these results in the finite abelian case, and the excellent survey by Folland and Sitaram~\cite{FollandSitaram}.

In this article, we study generalizations of the above bounds to general compact groups. An immediate difficulty in generalizing estimates such as~\eqref{finite-func} and~\eqref{finite-operator} to nonabelian groups is to settle upon an appropriate interpretation of $|\supp \hat{f}|$, as the Fourier transform is now a collection of linear operators.  We show that a natural analogue of (\ref{finite-operator}) still holds when $\widehat G$ is the collection of irreducible representations of a compact group $G$. We also refine this bound by establishing an uncertainty principle for a wider class of operators, that may operate inside the spaces of the various irreducibles representations of $G$. On finite groups, this principle subsumes the results of Donoho and Stark~\cite{Donoho} and Meshulam~\cite{Meshulam}. A corollary of our principle appropriately generalizes the statement~\eqref{finite-func} to compact groups, improving upon a result of Matolcsi and Sz\"ucs~\cite{Matolcsi}. This corollary is also an improvement, in the setting of compact groups, over the principle for which Kaniuth~\cite{Kaniuth} recently characterized all minimizing functions. Our primary contribution is the following theorem.

\begin{intro-thm}\label{thm:main}
Let $G$ be a compact group with Haar measure $\mu$, and let $\op P$ and $\op R$ be operators on $L^2(G)$. If $\;\op P$ commutes with projection onto every measurable subset of $G$ and $\op R$ commutes with left multiplication by elements of $G$, then $\|\op {PR}\|_2 = \|\op P \cdot \chi_G \|_2 \|\op R\|_2$.
\end{intro-thm}
\noindent Here $\chi_G$ denotes the characteristic function of $G$ and $\| \cdot \|_2$ denotes the $L^2$-norm of a function as well as the Hilbert-Schmidt norm of an operator. An immediate consequence of the above is the operator-theoretic uncertainty principle $\|\op{PR}\| \leq \|\op P \cdot \chi_G \|_2 \|\op R\|_2$, where $\| \cdot \|$ denotes operator norm. We remark that Theorem \ref{thm:main} still holds when ``left'' is replaced by ``right,'' corresponding to a different choice of Fourier transform. Theorem \ref{thm:main} allows us to prove the following uncertainty principle concerning 
functions on a compact group.
\begin{intro-thm}\label{thm:func}
Let $G$ be a compact group with Haar measure $\mu$, and $f$ a nonzero element of $L^2(G)$. Then
$\mu(\supp f) \cdot \sum_{\rho \in \hat G} d_\rho \rank \hat f(\rho) \geq 1$.
\end{intro-thm}
\noindent This fact was first proved for finite groups by Russell and Shparlinski~\cite{Shparlinski}, and improves upon the previously known result of Matolcsi and Sz\"ucs \cite{Matolcsi}, where the size of the support of $\hat f$ was given by $\sum_{\{\rho : \hat f(\rho) \neq 0\}} d_\rho^2$. We remark that there is a large family of examples where the principle in Theorem \ref{thm:func} is tight, but the one due to Matolcsi and Sz\"ucs is not; specifically, this is the case whenever $f$ is the characteristic function of a non-normal subgroup of finite index. Theorem \ref{thm:main} also implies the following result on finite groups.
\begin{intro-cor}\label{cor:cor1}
Let $\op P$ and $\op R$ be projection operators on the group algebra of a finite group $G$. 
If $\;\op P$ commutes with projection onto elements of $G$, and $\op R$ commutes with left-multiplication by elements of $G$, then $\|\op {PR}\|^2 \leq (\rank \op P \rank \op R)/|G|$.
\end{intro-cor}
\noindent This corollary, in turn, easily implies a result of the form (\ref{finite-operator}) for all finite groups.
\begin{intro-cor}\label{cor:cor2}
Let $G$ be a finite group, and let $S \subset G$ and $T \subset \hat G$. If $\op P_S$ denotes projection onto $S$ (in the group basis), and $\op R_T$ denotes projection onto $T$ (in the Fourier basis), then
$$\
\|\op {P_S R_T}\|^2 \leq \frac{|S| \cdot \sum_{\rho \in T} d_\rho^2}{|G|}~.
$$
\end{intro-cor}
\noindent In the next section, we discuss the so-called time-limiting and band-limiting operators appearing in the statements above. We then prove Theorems \ref{thm:main} and \ref{thm:func}, and discuss the implications in the setting of finite groups. We close with an alternative proof of Corollary \ref{cor:cor1}, using only basic results from the representation theory of finite groups.

\section{Preliminaries}

\subsection{Time-limiting operators}

The definitions of \emph{time-limiting} and \emph{band-limiting} operators, which act on the space of square-integrable functions on a group, are inspired by the signal processing applications of the uncertainty principle considered by Donoho and Stark~\cite{Donoho}. In what follows, $G$ will denote a compact group with Haar measure $\mu$, normalized so that $\mu G = 1$. A time-limiting projection operator $\op P_S$ acts on $L^2(G)$ simply by projecting to a subset $S$ of the group. More generally, we may consider the linear operator $\op P_f$ given by pointwise multiplication by a bounded function $f$. While the operator norm of $\op P_f$ is equal to the $L^\infty$-norm of $f$, the Hilbert-Schmidt norm of $\op P_f$ is unbounded unless $f = 0$ almost everywhere. The following proposition shows that these operators are a natural generalization of time-limiting projections.

\begin{prop}\label{prop: time-operator}
Let $G$ be a compact group. A bounded linear operator $\op P$ on $L^2(G)$ commutes with projection onto every measurable subset of $G$ if and only if $\op P = \op P_f$ for some bounded $f \in L^2(G)$.
\end{prop}
\begin{proof}
The reverse direction is clear. For the forward direction, we first define $f = \op P \cdot \chi_G$ and note that $\op P$ acts by pointwise multiplication by $f$ on all simple functions. This is enough to show that $f$ must also be bounded. Supposing otherwise, let $S_k = \{x : f(x) \in [k, k+1]\}$ and define a simple function $g$ by
$$
g(S_k) = \begin{cases} \frac{1}{k \sqrt{\mu(S_k)}} &\text{ if } \mu(S_k) > 0 \\
				0 &\text{ otherwise.}
	\end{cases}
$$
Evidently, $g \in L^2(G)$, but
$$
\|P \cdot g\|_2^2 = \int_G |f(x)g(x)|^2 d\mu(x) \geq \sum_k \mu(S_k) \cdot k |g(S_k)|^2~,
$$ 
and an infinite number of terms in this sum are equal to $1$. Hence $\op P \cdot g \notin L^2(G)$, and thus $f$ must be bounded. Finally, an arbitrary function $g \in L^2(G)$ is an $L^2$-limit of simple functions $g_n$. Since $\op P$ is bounded (and hence continuous), and $f$ is a bounded function, we have
$$
\op P \cdot g = \op P \cdot \lim_{n \rightarrow \infty} g_n = \lim_{n \rightarrow \infty} \op P \cdot g_n = \lim_{n \rightarrow \infty} fg_n = fg~.
$$
\end{proof}

\subsection{Band-limiting operators}
Recall that the Fourier transform of a function $f \in L^2(G)$ at the irreducible representation $\rho$
is given by
\begin{equation}\label{eq:fourier-transform}
\hat f(\rho) = \int f(x) \rho(x)^\dagger~d\mu(x)\,.
\end{equation}
Using the Fourier inversion formula, we may then define band-limiting projection onto a subset $T$ of $\hat G$ by setting
$$
\op R_T \cdot f(x) = \sum_{\rho \in T} d_\rho \trace \left[ \hat f(\rho) \rho(x) \right]\,.
$$
We remark that these operators commute with both the left and the right action of $G$ on $L^2(G)$. In terms of the irreducible decomposition of $L^2(G)$ given by the Peter-Weyl theorem, $\op R_T$ operates on $L^2(G)$ by projecting onto the subspace $\bigoplus_{\rho \in T} \mathcal E_\rho$, where $\mathcal E_\rho$ denotes the $\rho$-isotypic subspace. While this is satisfactory for studying abelian groups, where $\dim \mathcal E_\rho = 1$, such operators become increasingly coarse as the dimensions of the various irreducible representations of $G$ increase. For this reason, we wish to consider more general band-limiting operators in the non-abelian setting. A general band-limiting operator will be described by a collection $\op R = \{\op R_\rho\}_{\rho \in \hat G}$ of linear operators, where $\op R_\rho$ operates on the space of $\rho$. The action of such an operator on a function $f \in L^2(G)$ is then given by
\begin{equation}\label{eq:left-commuting}
\op R \cdot f(x) = \sum_{\rho \in \hat G} d_\rho \trace \left[ \op R_\rho \hat f(\rho) \rho(x) \right]~.
\end{equation}
The previous notion of a band-limiting operator $\op R_T$ is the special case where $R_\rho = \identity_\rho$
for $\rho \in T$ and $R_\rho = 0$ otherwise. 

Given a group element $x \in G$, let $\op L_x$ denote the operator corresponding to the left action of $x$ on $L^2(G)$, i.e., $[\op L_x\cdot f](y) = f(x^{-1}y)$. If an operator $\op R$ satisfies $\op {RL}_x = \op L_x\op R$ for every $x \in G$, then we say that $\op R$ commutes with left multiplication.

\begin{prop}\label{prop: operator-description}
Let $G$ be a compact group, and $\op R$ a linear operator on $L^2(G)$. Then $\op R$ commutes with left multiplication if and only if $\;\op R \cdot f(x) = \sum_{\rho \in \hat G} d_\rho \trace [\op R_\rho \hat f(\rho) \rho(x)]$ for some collection $\{\op R_\rho\}_{\rho \in \hat G}$, where each $\op R_\rho$ is a linear
operator on the space of $\rho$.
\end{prop}
\begin{proof}
We first observe that
$$
\widehat{L_y \cdot f} (\rho) = \int f(y^{-1}x) \rho(x)^\dagger d\mu(x) 
= \int f(x) \rho(yx)^\dagger d\mu(x) = \hat f (\rho) \rho(y^{-1})~.
$$
If a linear operator $\op R$ satisfies (\ref{eq:left-commuting}), then
$$
\left[\op {RL}_y \cdot f\right] (x)
= \sum_{\rho \in \hat G} d_\rho \trace \left[ \op R_\rho \hat f (\rho) \rho(yx) \right]
= \left[\op L_y\op R \cdot f\right](x),
$$
i.e. $\op R$ commutes with left multiplication.

Now let $\op R$ be any linear operator on $L^2(G)$ which commutes with left multiplication. By Schur's Lemma, such an operator decomposes into a direct sum of linear operators $\op A_\rho$, one for each irreducible representation $\rho \in \hat G$. Each $\op A_\rho$ acts on the $\rho$-isotypic subspace $\mathcal E_\rho$ from the decomposition of $L^2(G)$ given by the Peter-Weyl theorem. Schur's Lemma also asserts that each $\op A_\rho$ is determined by a linear operator $\op R_\rho$ acting on a single space of the representation $\rho$. The action of $\op A_\rho$ on the Fourier transform of $f$ is precisely by matrix multiplication by $\op R_\rho$ on the left. Taken together with the Fourier inversion formula, this means precisely that $\op R$ has the form (\ref{eq:left-commuting}).
\end{proof}

\section{Results}
%%%%%%%%%%%%%%%%%%%%%%%%%%%%%%%%%%%%%%%%%%%%%%%%%%%%%%%%%%%%%%%%%%%%%%%%%
%%%%%%%%%%%%%%%%%%%%%%%%%%%%%%%%%%%%%%%%%%%%%%%%%%%%%%%%%%%%%%%%%%%%%%%%%

\subsection{Operator uncertainty principles} Our primary technical contribution is the following.

\begin{thm} \label{thm1}
Let $G$ be a compact group, and $f$ a bounded measurable function on $G$. 
Let $\op R$ be an operator on $L^2(G)$ which commutes with right
multiplication by elements of $G$. Then
$$
\|\op P_f\op R\|_2^2 = \|f\|_2^2 \sum_{\rho \in \hat G} d_\rho \|\op R_\rho\|_2^2~,
$$
where $R_\rho$ is the operator on the representation space of $\rho$ from the decomposition of $R$ implied by Proposition \ref{prop: operator-description}.
\end{thm}
\begin{proof}
For every irreducible representation $\rho$ of $G$, let $\Pi_\rho$ denote projection onto the $\rho$-isotypic subspace $\mathcal E_\rho \subset L^2(G)$. By the Peter-Weyl theorem, the $\Pi_\rho$ form a resolution $\sum_\rho \Pi_\rho = \identity$ of the identity operator on $L^2(G)$. We can thus express the Hilbert-Schmidt norm of $\op P_f \op R$ as follows:
\begin{align*}
\|\op P_f\op R\|_2^2 & = \trace \left[\op R^\dagger \op P_f^\dagger \op P_f \op R \right]
		= \trace \left[\sum_\rho \Pi_\rho \op R^\dagger \op P_{|f|^2} \op R \sum_\sigma \Pi_\sigma \right] \\
		& = \sum_\rho \sum_\sigma \trace \left[\Pi_\rho \op R^\dagger \op P_{|f|^2} \op R \Pi_\sigma \right]
		= \sum_\rho \trace \left[\Pi_\rho \op R^\dagger \op P_{|f|^2} \op R \Pi_\rho \right]~,
\end{align*}
where we have made use of the fact that $(\op P_f)^\dagger = \op P_{\bar f}$. We remark that the operators $\Pi_\rho$ also commute with left multiplication by elements of $G$. Indeed, their images are invariant subspaces of $L^2(G)$ viewed as a representation of $G$ under the left-multiplication action. We thus have, for every $x \in G$,
\begin{align*}
\|\op P_f\op R\|_2^2 & = \sum_\rho \trace \left[\op L_x\op L_{x^{-1}}\Pi_\rho \op R^\dagger \op P_{|f|^2} \op R \Pi_\rho \right] \\
  & = \sum_\rho \trace \left[\op L_{x^{-1}}\Pi_\rho \op R^\dagger \op P_{|f|^2} \op R \Pi_\rho \op L_x \right]  \\
	& = \sum_\rho \trace \left[\Pi_\rho \op R^\dagger \op L_{x^{-1}}\op P_{|f|^2}\op L_x \op R \Pi_\rho \right]~.
\end{align*}
Integrating both sides over $G$, we conclude that
\begin{align*}
\|\op P_f\op R\|_2^2 & = \int_x \|\op P_f\op R\|_2^2 ~d\mu(x) \\
		& = \int_x \sum_\rho \trace \left[\Pi_\rho \op R^\dagger \op L_{x^{-1}} \op P_{|f|^2}\op L_x \op R \Pi_\rho \right] ~d\mu(x) \\
		& = \sum_\rho \trace \left[\Pi_\rho \op R^\dagger \int_x \op L_{x^{-1}} \op P_{|f|^2}\op L_x ~d\mu(x) \op R \Pi_\rho \right]~.
\end{align*}
The last integral above denotes the so-called \emph{weak operator integral}. One may arrive at the proper definition of such an integral, for instance, by requiring that the integral $\int \langle \Phi(x)u, v \rangle$ of each matrix element is equal to the corresponding matrix element of the integral $\int \Phi(x)$. This is well-defined for any bounded operator-valued function $\Phi$, i.e., if there exists $M$ such that $\|\Phi(x)\|\leq M$ for every $x$. Moreover, the integral thus defined commutes with taking the trace, and commutes with composition with bounded operators on the left and the right; in fact, this is precisely the integral appearing in the operator-valued Fourier transform \eqref{eq:fourier-transform} for compact groups. For further details on issues of integrability and measurability of operator-valued functions, consult Fell and Doran~\cite{FellDoran} or Conway~\cite{Conway}. Returning to our calculation, we now wish to show that the operator $\int_G \op L_{x^{-1}} \op P_{|f|^2} \op L_x~d\mu(x)$ is in fact just scaling by $\|f\|_2^2$, i.e., the operator $\op P_{\|f\|_2^2}$. Indeed, we see that for any $g, h \in L^2(G)$,
\begin{align*}
\left \langle \int_x \op L_{x^{-1}} \op P_{|f|^2} \op L_x~d\mu(x) \cdot g, ~h \right \rangle
& = \int_x \left\langle \op L_{x^{-1}} \op P_{|f|^2} \op L_x \cdot g, ~h \right \rangle ~d\mu(x) \\
& = \int_x \int_y \op L_{x^{-1}} \op P_{|f|^2} \op L_x \cdot g(y)\overline{h(y)} ~d\mu(y)~d\mu(x) \\
& = \int_x \int_y |f(xy)|^2 g(y)\overline{h(y)} ~d\mu(y)~d\mu(x) \\
& = \int_y \int_x |f(xy)|^2 d\mu(x) g(y)\overline{h(y)} ~d\mu(y) \\
& = \left \langle \|f\|_2^2 \cdot g, h \right \rangle~,
\end{align*}
where the second-to-last step follows from Fubini's Theorem. Finally, we have
\begin{align*}
\|\op P_f\op R\|_2^2 & = \sum_\rho \trace \left[\Pi_\rho \op R^\dagger \op P_{\|f\|_2^2} \op R \Pi_\rho \right] 
		 = \|f\|_2^2 \sum_\rho \trace \left[\Pi_\rho \op R^\dagger \op R \Pi_\rho \right] \\
		& = \|f\|_2^2 \sum_\rho \|\op R \Pi_\rho\|_2^2
		 = \|f\|_2^2 \sum_\rho \|\op A_\rho\|_2^2
		 = \|f\|_2^2 \sum_\rho d_\rho \|\op R_\rho\|_2^2~,
\end{align*}
where the last two steps follow from the proof of Proposition \ref{prop: operator-description}~. 
\end{proof}

\noindent The following operator uncertainty principle follows easily from the above theorem.
\begin{thm-restate-main}
Let $G$ be a compact group, and let $\op P$ and $\op R$ be operators on $L^2(G)$~. If
$\;\op P$ commutes with projection onto every measurable subset of $G$ and
$\op R$ commutes with left multiplication by elements of $G$, then
$$
\|\op {PR}\| \leq \|\op P \cdot \chi_G \|_2 \|\op R\|_2~.
$$
\end{thm-restate-main}
\begin{proof}
By Proposition \ref{prop: time-operator}, $\op P = \op P_f$, where $f = \op P \cdot \chi_G \in L^\infty(G)$.
Theorem \ref{thm1} then asserts that
$$
\|\op {PR}\|_2^2 = \|\op P \cdot \chi_G \|_2^2 \sum_\rho d_\rho \|\op R_\rho\|_2^2 = \|\op P \cdot \chi_G \|_2^2 \|\op R\|_2^2~.
$$
The result now follows from the fact that the operator norm is bounded above by the Hilbert-Schmidt norm.
\end{proof}

\subsection{Function uncertainty principles}

In this section, we apply Theorem \ref{thm:main} to prove a ``classical'' uncertainty principle along the lines of (\ref{finite-func}), relating the support of a nonzero function on a compact group  to the support of its Fourier transform. Since we are given a normalized measure $\mu$ on $G$, there is a canonical and natural way to quantify the size of the support of $f$, i.e., $|\supp f| = \mu(\supp f)$. The size of the support of $\hat f$, on the other hand, involves a nontrivial choice of dual measure on $\hat G$. For finite groups, a natural choice is \emph{Plancherel measure}, which assigns mass $d_\rho^2/|G|$ to each representation $\rho$. As this is the dimensionwise fraction of the group algebra consisting of irreps isomorphic to $\rho$, the Plancherel  measure of $\hat G$ is equal to the (normalized) Haar measure of $G$. For general compact groups, we will assign measure $d_\rho^2$ to each representation $\rho$. This is still a natural choice, since $d_\rho^2$ is equal to the dimension of the $\rho$-isotypic subspace of $L^2(G)$. We thus set, for $T \subset \hat G$, $|T| = \sum_{\rho \in T} d_\rho^2$; meanwhile, for $S \subset G$, we set $|S| = \mu S$. Now, consider the following statement:
\begin{equation}\label{QUP}
\text{\emph{If~}} |\supp f| < |G| \text{\emph{ and }} |\supp \hat f| < |\hat G|\text{\emph{, then }} f = 0\,.
\end{equation}
We emphasize that for general compact groups, and indeed in all cases of interest for (\ref{QUP}), the quantity $|\hat G|$ is infinite. Hogan~\cite{Hogan} showed that if $G$ is infinite and compact, then (\ref{QUP})  is valid if and only if $G$ is connected. Echterhoff, Kaniuth and Kumar~\cite{EchterhoffKaniuthKumar} showed that if $G$ has a noncompact, nondiscrete, cocompact normal subgroup $H$ that satisfies (\ref{QUP}), then $G$ satisfies it as well. Another principle, given in Kutyniok \cite{Kutyniok}, states that
\begin{equation}\label{weakQUP}
\text{\emph{If~}} |\supp f||\supp \hat f|_1 < 1\text{\emph{, then }} f = 0~,
\end{equation}
where $|\supp \hat f|_1 = \sum_{\rho \in \supp \hat f} d_\rho$. Kutyniok proves that a compact group 
$G$ satisfies (\ref{weakQUP}) if and only if $G$ modulo the connected component of the identity is an abelian group.
While the choice between $|\cdot|$ and $|\cdot|_1$ has no impact on abelian groups, or on the principle
(\ref{QUP}), it is quite significant for the principle (\ref{weakQUP}) on arbitrary compact groups.
For instance, if $G$ has a finite quotient $G/H$, then $|\supp \chi_H| = |H|$ while
$$
|\supp \hat \chi_H| = \sum_{\rho \in \supp \hat \chi_H} d_\rho^2 = \sum_{\rho \in \widehat{G/H}} d_\rho^2 = [G:H] = 1/|H|.
$$
Hence $|\supp \chi_H||\supp \hat \chi_H| = 1$, while any smaller choice of dual measure would result in inequality. In fact, the entire proof of the main theorem of \cite{Kutyniok} still holds if we assign measure $1$ to one-dimensional irreps, and measure $0$ to any $\rho$ satisfying $d_\rho > 1$. 

For arbitrary compact groups, our Theorem \ref{thm1} implies the following result.
\begin{thm-restate-func}
Let $G$ be a compact group, and $f$ a nonzero element of $L^2(G)$. Then
$$
|\supp f| \sum_{\rho \in \hat G} d_\rho \rank \hat f(\rho) \geq 1~.
$$
\end{thm-restate-func}
\begin{proof}
Let $f$ be a nonzero element of $L^2(G)$. For each $\rho \in \hat G$, let $\op R_\rho$ be the operator on the space of $\rho$ which projects to ${\bf Im} \hat f(\rho)$, as a subspace of the representation space of $\rho$. Define the operator $\op R$ on $L^2(G)$ by
$$
\op R \cdot g(x) = \sum_\rho d_\rho \trace \left[ \op R_\rho \hat g(\rho) \rho(x) \right]
$$
and recall that $\op R$ commutes with left multiplication by Proposition \ref{prop: operator-description}. Set $\op P = \op P_{\supp f}$, i.e., $L^2$ projection onto the support of $f$. By Theorem \ref{thm1},
$$
\|\op{PR}\|_2^2 = |\supp f| \sum_{\rho \in \hat G} d_\rho \|\op R_\rho\|_2^2
= |\supp f| \sum_{\rho \in \hat G} d_\rho \rank \hat f(\rho)\,.
$$
As both $\op P$ and $\op R$ are orthogonal projections, and $\op{PR} \cdot f = f$, we have that $\|\op{PR}\| = 1$. The claim follows from the fact that the operator norm is bounded above by the Hilbert-Schmidt norm.
\end{proof}
\noindent We now provide a simple proof of the above in the special case where $f$ is a bounded function and achieves a global maximum somewhere on the group. The proof is a straightforward adaptation to the compact case of the proof due to Russell and Shparlinski ~\cite{Shparlinski}.
\begin{proof}
Let $s \in G$ be the group element where $f$ achieves its maximum. By the inversion formula,
$$
\|f\|_2^2 \leq |\supp f|\|f\|_\infty^2 = |\supp f|\left|\sum_{\rho \in \hat G} d_\rho \trace \left[ \hat f(\rho) \rho(s) \right] \right|^2 = |\supp f|\left|\sum_{\rho \in \hat G;~1 \leq i,j \leq d_\rho} \left(\hat f(\rho) \rho(s)\right)_{jj} \right|^2
$$
Noting that the sum above has $\sum_{\rho \in \hat G} d_\rho \rank \hat f(\rho)$ many terms, we apply Cauchy-Schwarz:
\begin{align*}
\|f\|_2^2 
&\leq |\supp f|\left(\sum_{\rho \in \hat G} d_\rho \rank \hat f(\rho)\right) \sum_{\rho \in \hat G;~1 \leq i,j \leq d_\rho} \left|\left(\hat f(\rho) \rho(s)\right)_{jj}\right|^2 \\
&\leq |\supp f|\left(\sum_{\rho \in \hat G} d_\rho \rank \hat f(\rho)\right) \sum_{\rho \in \hat G} d_\rho \trace\left[\hat f(\rho) \rho(s) \rho(s)^\dagger \hat f(\rho)^\dagger\right]\\
&= |\supp f|\left(\sum_{\rho \in \hat G} d_\rho \rank \hat f(\rho)\right) \sum_{\rho \in \hat G} d_\rho \|\hat f(\rho)\|_2^2~.
\end{align*}
The claim now follows from the Plancherel identity.
\end{proof}
\noindent An obvious consequence of the above theorem is that any nonzero $f \in L^2(G)$ satisfies
\begin{equation}\label{oldup}
|\supp f|\sum_{\rho \in \supp \hat f} d_\rho^2 \geq 1~.
\end{equation}
Indeed, this is precisely the uncertainty principle appearing in~\cite{Matolcsi}, and the best previously known multiplicative uncertainty principle for functions on arbitrary compact groups. 

We now produce a family of natural examples where our principle from Theorem \ref{thm:func} is tight, but \eqref{oldup} is not. Let $H$ be a subgroup of $G$ with finite index, and note that the scaled Fourier transform $\widehat{\chi_H}/\mu H$ of the characteristic function of $H$ is a projection operator. By Plancherel, it follows that our principle is tight for these functions:
$$
\mu(\supp \chi_H) \sum_{\rho \in \hat G} d_\rho \rank \widehat{\chi_H} (\rho)
= \mu H \sum_{\rho \in \hat G} d_\rho \rank\left[\frac{\widehat{\chi_H} (\rho)}{\mu H}\right]
= \frac{1}{\mu H} \sum_{\rho \in \hat G} d_\rho \|\widehat{\chi_H}(\rho)\|_2^2
= 1~.
$$
As we now show, the principle \eqref{oldup} is not tight for $\chi_H$ unless $H$ is normal.
\begin{prop}
Let $H$ be a subgroup of a compact group $G$ with $[G:H]$ finite. Then $H$ is normal in $G$ if and only if, for every $\rho \in \hat G$, $\rank \widehat{\chi_H}(\rho)$ is either $d_\rho$ or zero.
\end{prop}
\begin{proof}
Assume first that $H$ is normal, and let $H^\perp = \{\rho \in \hat G : H \subset \ker \rho\}$. It's clear that $\widehat{\chi_H}(\rho) = \identity_\rho$ if $\rho \in H^\perp$. Summing over $H^\perp \cong \widehat{G/H}$, we have
$$
\sum_{\rho \in H^\perp} d_\rho \|\widehat{\chi_H}(\rho)\|_2^2 = \sum_{\rho \in \widehat{G/H}} (\mu H)^2 d_\rho^2 = |G/H| (\mu H)^2 = \mu G \mu H = \mu H~;
$$
the Plancherel equality then implies that $H^\perp$ is in fact the entire support of $\widehat{\chi_H}$.

Now let us suppose that the ranks of the Fourier transforms of $\chi_H$ are all either full or zero. We claim that $H$ is equal to the intersection of the kernels of the irreps for which $\rank \widehat{\chi_H}(\rho) = d_\rho$. We will denote the set of all such irreps by $H^\perp$. Given any $h \in H$ and $\rho \in H^\perp$, we have $\rho(h) \widehat{\chi_H}(\rho) = \widehat{\chi_H}(\rho)$; together with the fact that $\widehat{\chi_H}(\rho)$ has full rank, this implies that $\rho(h) = \identity_\rho$. For the other containment, we apply the Fourier inversion formula to compute the value of $\chi_H$ at a group element $g \in \bigcap_{\rho \in H^\perp} \ker \rho$~:
$$
\chi_H(g) = \sum_{\rho \in \hat G} d_\rho \trace \left[ \widehat{\chi_H}(\rho) \rho(g)\right]
= \frac{1}{\mu H} \sum_{\rho \in H^\perp} d_\rho \trace \left[\identity_\rho \rho(g)\right]
= \frac{1}{\mu H} \sum_{\rho \in H^\perp} d_\rho^2
= 1~,
$$
which completes the proof.
\end{proof}

\subsection{Consequences in finite groups}

In this section, we prove a corollary of the main theorem in the setting of finite groups.

\begin{cor-restate-cor1}
Let $\op P$ and $\op R$ be projection operators on the group algebra of a finite group $G$. If $\;\op P$ commutes with projection onto elements of $G$, and $\op R$ commutes with left-multiplication by elements of $G$, then
$\|\op {PR}\|^2 \leq \rank \op P \rank \op R / |G|~.$
\end{cor-restate-cor1}
\begin{proof}
A projection operator on $\CC G$ commutes with projection onto elements of $G$ if and only if it projects onto some set $S \subset G$, in which case its rank is $|S|$. Meanwhile, a projection operator commutes with left-multiplication exactly when it respects the decomposition of $\CC G = L^2(G)$ into irreducible spaces according to the left action. Since the restrictions $\op R_\rho$ of $\op R$ to these irreducible spaces are themselves projection operators, we have
$$
\rank \op R = \sum_{\rho \in \hat G} d_\rho \cdot \rank \op R_\rho = \sum_{\rho \in \hat G} d_\rho \cdot \| \op R_\rho\|_2^2~.
$$
The result now follows from Theorem \ref{thm1}.
\end{proof}

\noindent We now provide an alternative proof of the above, making use of only basic results from the representation theory of finite groups.

\begin{proof}
We first claim that
\begin{equation}\label{finite-claim}
\sum_{x \in G} {\bf tr}\left(\op R_{\rho_1}\rho_1(x)\right)\overline{{\bf tr}\left(\op R_{\rho_2}\rho_2(x)\right)} = 
\begin{cases}
0 & \text{if } \rho_1 \neq \rho_2 \\
\frac{{\bf rk}\op R_{\rho_1} |G|}{d_{\rho_1}} & \text{if } \rho_1 = \rho_2.
\end{cases}
\end{equation}
When $\rho_1 \neq \rho_2$, the result follows from the orthogonality of the Fourier basis functions. The two trace terms are linear combinations of such functions, and hence their product is zero. Now suppose the two irreps are equal, and let $\rho = \rho_1 = \rho_2$. Then
$$
\sum_{x \in G} {\bf tr}\left(\op R_\rho \rho(x)\right)\overline{{\bf tr}\left(\op R_\rho \rho(x)\right)} 
= \text{ {\bf tr}}\left(\op R_\rho \otimes \op R_\rho^\dagger \cdot \sum_{x \in G} \rho(x) \otimes \rho(x^{-1})\right).
$$
We observe the following property of the operator $\op J = \sum_{x \in G} \rho(x) \otimes \rho(x^{-1}):$
$$
\rho(h)\otimes\rho(g) \cdot \op J  = \sum_{x \in G} \rho(hx) \otimes \rho(gx^{-1})
 = \sum_{y \in G} \rho(yg) \otimes \rho(y^{-1}h) = \op J \cdot \rho(g)\otimes\rho(h)\enspace,
$$
where we have used the substitution $y = hxg^{-1}.$ Letting $\op K$ be the interchange operator $\op K:a \otimes b \mapsto b \otimes a$, we have $\rho(h)\otimes\rho(g) \cdot \op{JK} = \op{JK} \cdot \rho(h)\otimes\rho(g)$, which by Schur's Lemma means that $\op{JK}$ is a homothety. Since $\op{K}$ is its own inverse, we conclude that $\op{J} = \lambda \op{K}$ for some scalar $\lambda$. It's easy to check that $\trace \op J = |G|$ while $\trace \op K = d_\rho$, so that $\lambda = |G|/d_\rho.$ Now let $B$ be an orthonormal basis for the space of $\rho$ consistent with the projector $\op R_\rho$. Returning to our original computation, we now have
\begin{align*}
 \text{ {\bf tr}}\left(\op R_\rho \otimes \op R_\rho^\dagger \cdot \op J \right) 
 & = \frac{|G|}{d_\rho} \sum_{a, b \in B} \bra a \otimes b, \op R_\rho b \otimes \op R_\rho a \ket
   = \frac{|G|}{d_\rho} \sum_{a, b \in B} \bra a, \op R_\rho b \ket \bra b, \op R_\rho a \ket \\
 & = \frac{|G|}{d_\rho} \sum_{a \in {\bf im} \op R_\rho} \bra a, a \ket \bra a, a \ket
   = \frac{|G|\rank \op R_\rho}{d_\rho}~,
\end{align*}
establishing our claim. We now compute the matrix entries of $\op{RP}$ in the basis of delta functions on $G$, recalling that $P$ projects to some subset $S \subset G$.
\begin{align*}
\op{RP} \cdot f(x) 
&= \sum_{\rho \in \hat G} d_\rho \trace \left(\op R_\rho\rho(x)^\dagger \widehat {P\cdot f}(\rho)\right) 
 = \sum_{\rho \in \hat G} \frac{d_\rho}{|G|} \trace \left(\op R_\rho\rho(x)^\dagger \sum_{y \in S} f(y) \rho(y)\right)\\
& = \sum_{\rho \in \hat G} \frac{d_\rho}{|G|}\sum_{y \in S} \trace \left(\op R_\rho\rho(x)^\dagger \rho(y)\right)f(y) 
 = \sum_{y \in G} \sum_{\rho \in \hat G} \frac{d_\rho\chi_S(y)}{|G|} \trace \left(\op R_\rho\rho(x^{-1}y)\right) f(y)~.
\end{align*}
We then compute the $L^2$-norm directly:
\begin{align*}
\|\op{RP}\|_2^2 &= \sum_{x,y \in G}\left|\sum_{\rho \in \hat G} \frac{d_\rho\chi_S(y)}{|G|} \textbf{tr}\left(\op R_\rho\rho(x^{-1}y)\right)\right|^2 \\
&= \frac{1}{|G|^2}\sum_{x,y \in G} \sum_{\rho_1, \rho_2 \in \hat G} d_{\rho_1} d_{\rho_2}\chi_S(y) 
                   \trace \left(\op R_{\rho_1}\rho_1(x^{-1}y)\right)\overline{\trace \left(\op R_{\rho_2}\rho_2(x^{-1}y)\right)} \\
&= \frac{1}{|G|^2}\sum_{y \in S} \sum_{\rho_1, \rho_2 \in \hat G} d_{\rho_1} d_{\rho_2} \sum_{x \in G}
                   \textbf{tr}\left(\op R_{\rho_1}\rho_1(x^{-1}y)\right)\overline{\textbf{tr}\left(\op R_{\rho_2}\rho_2(x^{-1}y)\right)}\enspace.
\end{align*}
By \eqref{finite-claim}, the innermost sum is zero except when $\rho_1 = \rho_2$, yielding
$$
\|\op{RP}\|_2^2 = \frac{1}{|G|^2}\sum_{y \in S} \sum_{\rho \in \hat G} d_\rho^2 \frac{\rank \op R_\rho |G|}{d_\rho} 
 = \frac{|S| \sum_{\rho \in \hat G} d_\rho\rank \op R_\rho}{|G|} = \frac{\rank \op R \rank \op P}{|G|}~.
$$
\end{proof}

\noindent We conclude by remarking that if $\op P$ projects to $S \subset G$ and $\op R$ projects to $T \subset \hat G$, then the above proves that $\|\op {PR}\|^2 \leq (\sum_{\rho \in T} d_\rho^2)|S|/|G|$, which is the nonabelian analogue of the original statement \eqref{finite-operator} of Donoho and Stark~\cite{Donoho}.

\end{document}